\documentclass[a4paper,12pt]{article}
\usepackage[english]{babel}
\usepackage[T1]{fontenc}
\usepackage[utf8]{inputenc}
\usepackage{amsthm}
\usepackage[tbtags]{amsmath}
\usepackage{amsfonts,amssymb}
\usepackage[usenames]{color}
\usepackage[normalem]{ulem}
\begin{document}

\newtheorem{theorem}{Theorem}
\newtheorem{lemma}{Lemma}
\newtheorem{proposition}{Proposition}
\newtheorem{Cor}{Corollary}

\begin{center}
\textbf{Minimal Finite Non-Commutative Centrally Essential Rings}
\end{center}

\begin{center}
A.N. Abyzov\footnote{Kazan (Volga Region) Federal University; email: adel.abyzov@kpfu.ru},
O.V. Lyubimtsev\footnote{Lobachevsky State University of Nizhny Novgorod; email: oleg\_lyubimcev@mail.ru},
A.A. Tuganbaev\footnote{National Research University MPEI, Lomonosov Moscow State University; email: tuganbaev@gmail.com}
\end{center}

\textbf{Abstract.} A ring $R$ with nonzero unity is called centrally essential if for every nonzero element $a$ of $R$ there exist nonzero central elements $x$ and $y$ such that $ax = y$. In this paper we obtain a complete classification, up to isomorphism, of non-commutative centrally essential rings of minimal order $64$.

\textbf{Keywords:} centrally essential ring, finite local ring.

The work of Adel Abyzov is supported by the Russian Science Foundation (grant 23-21-10086) and performed under the development program of the Volga Region Mathematical Center (agreement 075-02-2025-1725/1).
The work of Oleg Lyubimtsev is supported by the Ministry of Education and Science of the Russian Federation, project H-044-0\_2026-2028. The study of Askar Tuganbaev is supported by the Russian Science Foundation, project 26-11-00001.

\textbf{MSC2020 subject classification 16D25, 16R99}

\section{Introduction}

All rings considered in this paper are associative with nonzero unity. A ring $R$ is called \textbf{centrally essential} if for every nonzero element $a$ of $R$ there exist nonzero central elements $x$ and $y$ such that $ax = y$.\footnote{It is clear that a ring $R$ with center $Z$ is centrally essential if and only if the module $R_{Z}$ is an essential extension of the module $Z_{Z}$.} The class of centrally essential rings is a strict extension of the class of commutative rings. In \cite[Proposition 3.3]{MT18} it is proved that semiprime centrally essential rings are commutative. In \cite[Example 5.2]{LT25IJ} a non-commutative centrally essential ring of order $128$ is constructed. This gave rise to the natural problem of the minimum possible order of a non-commutative centrally essential ring, which was solved by the authors of the present paper in \cite{AbLyubTug2026}: it was proved that this minimum order equals $64$, the basic structural invariants of rings of order $64$ were established, and two families of non-commutative centrally essential $\mathbb{F}_2$-algebras of order $64$ were described under the additional assumption $J(R)^2 \subseteq \mathrm{Ann}(J(R)) \cap Z(R)$. Moreover, in \cite[Remarks 4.1 and 4.2]{AbLyubTug2026} open questions were formulated concerning the classification, up to isomorphism, within the first family of $\mathbb{F}_2$-algebras and concerning the description of rings of order $64$ of characteristic $4$. The goal of the present paper is to obtain a complete classification of all non-commutative centrally essential rings of order $64$ up to isomorphism, without any a priori restrictions on the square of the radical or on the characteristic of the ring.

We denote by $Z(R)$ and $J(R)$ the center and the Jacobson radical of a ring $R$, respectively. It is well known that the Jacobson radical of a right Artinian ring is nilpotent. The least $n$ such that $J(R)^n = 0$ is called the \textbf{nilpotency index} of $J(R)$. We write $\mathrm{Ann}(S) = \{r\in R \mid rS = 0\}$ for the left annihilator of a set $S$ in the ring $R$. A ring $R$ is called \textbf{local} if $R/J(R)$ is a division ring.

For completeness of the exposition, we recall the results of \cite{AbLyubTug2026} establishing the minimality of order $64$ and the basic numerical parameters of such rings.

\textbf{1.1. Remark.} \cite[Lemma 2.1]{AbLyubTug2026} There are no finite non-commutative centrally essential rings of order $\le 32$.

\textbf{1.2. Remark.} \cite[Lemma 2.2]{AbLyubTug2026} There are no finite non-commutative centrally essential rings whose order is strictly greater than $32$ and strictly less than $64$.

\textbf{1.3. Remark.} \cite[Theorems 1.1 and 1.2]{AbLyubTug2026} The minimum order of a finite non-commutative centrally essential ring is $64$. If $R$ is a finite non-commutative centrally essential ring of order $64$, then:
\begin{enumerate}
    \item[{\bf 1.}] The ring $R$ is local.
    \item[{\bf 2.}] The quotient ring $R/J(R)$ is isomorphic to the field $\mathbb{F}_2$, and the Jacobson radical $J(R)$ has order $32$.
    \item[{\bf 3.}] The center $Z(R)$ has order $16$, the intersection $Z_0 = Z(R) \cap J(R)$ has order $8$, and the subgroup $W_0 = \operatorname{Ann}(J(R)) \cap Z_0$ has order $2$ or $4$.
\end{enumerate}

The main result of the present paper is Theorem 1.4.

\textbf{1.4. Theorem.}
Let $R$ be a non-commutative centrally essential ring of order $64$, $J = J(R)$, $Z = Z(R)$, $Z_0 = Z \cap J$, and $W_0 = \mathrm{Ann}(J) \cap Z_0$. Then:
\begin{enumerate}
    \item[{\bf 1.}] The relations $J^2 \subseteq Z_0$, $[J, J] \subseteq W_0$, $JZ_0 \subseteq W_0$, $J^3 \subseteq W_0$, $J^4 = 0$, and $Z_0^2 = 0$ hold. In particular, the nilpotency index of the radical $J$ equals $3$ (when $J^2 \subseteq W_0$) or $4$ (when $J^2 \nsubseteq W_0$).
    \item[{\bf 2.}] The identities $2Z_0 = 0$ and $2J^2 = 0$ hold. In particular, $Z_0 \cong \mathbb{F}_2^3$ is a three-dimensional vector space over the field $\mathbb{F}_2$, and the characteristic $\mathrm{char}(R)$ of the ring $R$ equals $2$ or $4$.
    \item[{\bf 3.}] If $\dim_{\mathbb{F}_2} W_0 = 2$ ($|W_0| = 4$), then $J^2 \subseteq W_0$ (and hence $J^3 = 0$).
\end{enumerate}

On the basis of the above results, a complete classification of all non-commutative centrally essential rings $R$ of order $64$ up to isomorphism is obtained; see Propositions 3.2, 3.4, 3.6 in the case $\mathrm{char}(R) = 2$ and Propositions 4.1, 4.2, 4.3 in the case $\mathrm{char}(R) = 4$.

\section{Proof of Theorem 1.4}

Throughout this section, $R$ is a non-commutative centrally essential ring of order $64$, $J = J(R)$, $Z = Z(R)$, $Z_0 = Z \cap J$, and $W_0 = \mathrm{Ann}(J) \cap Z_0$.

\textbf{2.1. Lemma.}
The additive quotient group $J/Z_0$ is an elementary abelian $2$-group of rank $2$: $J/Z_0 \cong C_2 \times C_2$. In particular, for every element $r \in J$ the inclusion $2r \in Z_0$ holds.

\textbf{Proof.} By Remark 1.3 we have $|J| = 32$ and $|Z_0| = 8$, hence $|J/Z_0| = 4$. If the group $J/Z_0$ were cyclic with generator $a + Z_0$, then every element of $J$ would be expressible as $ka + z$, where $k \in \mathbb{Z}$ and $z \in Z_0 \subseteq Z$. Then for any two elements $k_1 a + z_1, k_2 a + z_2 \in J$ the equality $[k_1 a + z_1, k_2 a + z_2] = k_1 k_2 [a, a] = 0$ would hold. Since $R = J \cup (1 + J)$, the commutativity of the radical $J$ would imply the commutativity of the whole ring $R$, which contradicts the hypothesis. Therefore $J/Z_0 \cong C_2 \times C_2$ and $2J \subseteq Z_0$. $\square$

\textbf{2.2. Lemma.}
The subgroup $Z_0$ is a two-sided central ideal in $R$, that is, $JZ_0 \subseteq Z_0$. Moreover, a local ring $R$ with $R/J \cong \mathbb{F}_2$ and $JZ_0 \subseteq Z_0$ is centrally essential if and only if
$$
\mathrm{Ann}_J(Z_0) = \{r \in J \mid rZ_0 = 0\} \subseteq Z_0.
$$
\textbf{Proof.} Let $r \in J \setminus Z_0$ and $z \in Z_0$. Since $z \in Z$, we have $rz = zr \in J$. Suppose that $rz \notin Z_0$. We first verify that $rz \not\equiv r \pmod{Z_0}$. Indeed, if $rz - r = z' \in Z_0$, then $r(z - 1) = z'$. Since $z \in J$, the element $z - 1$ is invertible in the local ring $R$, and $(z - 1)^{-1} \in Z$ because $z \in Z$. Then $r = z'(z - 1)^{-1} \in Z \cap J = Z_0$, which contradicts the choice $r \notin Z_0$. Therefore the cosets $r + Z_0$ and $rz + Z_0$ are two distinct nonzero elements of the group $J/Z_0 \cong C_2 \times C_2$ and hence generate $J/Z_0$. Thus the additive group $J$ is generated by the subgroup $Z_0$ and the elements $r, rz$. But since $z \in Z$, the elements $r$ and $rz$ commute: $r(rz) = r^2 z = (rz)r$. Since $Z_0 \subseteq Z$, it follows that the ideal $J = \langle Z_0, r, rz \rangle$ is commutative, and hence the ring $R = \langle 1, J \rangle$ is commutative --- a contradiction. Thus $JZ_0 \subseteq Z_0$.

Let us prove the second part of the statement. For any invertible element $r = 1 + j \notin J$ ($j \in J$) and any nonzero element $w \in W_0 \setminus \{0\}$ (which exists by Remark 1.3), we have $rw = (1 + j)w = w \in Z \setminus \{0\}$. For any nonzero central element $r \in Z_0 \setminus \{0\}$ we have $r \cdot 1 = r \in Z \setminus \{0\}$. Finally, for a noncentral element $r \in J \setminus Z_0$ the equality $0 \neq rx \in Z$ ($x \in Z$) is impossible when $x \notin J$ (for then $x \in R^* \cap Z$ and $r = (rx)x^{-1} \in Z \cap J = Z_0$), while for $x \in Z_0$, by the inclusion $JZ_0 \subseteq Z_0$ proved above, the condition $rx \in Z$ holds automatically. Therefore the existence of $x \in Z$ with $0 \neq rx \in Z$ for all $r \in J \setminus Z_0$ is equivalent to $rZ_0 \neq 0$ for all $r \in J \setminus Z_0$, that is, to $\mathrm{Ann}_J(Z_0) \subseteq Z_0$. $\square$

\textbf{2.3. Lemma.}
In the ring $R$ the inclusions $J^2 \subseteq Z_0$ and $[J, J] \subseteq W_0$ hold.

\textbf{Proof.} By Lemmas 2.1 and 2.2, $Z_0$ is an ideal in $R$ and $2J \subseteq Z_0$, so the quotient ring $\overline{J} = J/Z_0$ is a two-dimensional nilpotent algebra over the field $\mathbb{F}_2$. If in a two-dimensional nilpotent algebra $\overline{J}$ over a field the square is nonzero ($\overline{J}^2 \neq 0$), then $\dim_{\mathbb{F}_2} \overline{J}^2 = 1$ and there exists an element $\overline{a} \in \overline{J}$ such that $\{\overline{a}, \overline{a}^2\}$ is a basis of $\overline{J}$. Lifting $\overline{a}$ to an element $a \in J$, we obtain that the additive group $J$ is generated by the central ideal $Z_0$ and the commuting elements $a$, $a^2$. Then $J$ and $R$ are commutative --- a contradiction. Therefore $\overline{J}^2 = 0$, that is, $J^2 \subseteq Z_0$.

Let us now prove that $[x, y] \in W_0$ for any $x, y \in J$. Since $J^2 \subseteq Z_0$, we have $[x, y] = xy - yx \in Z_0$, so it suffices to verify that $[x, y]J = 0$. For any $z \in Z_0$, in view of the centrality of $z \in Z$ and the inclusion $xz \in JZ_0 \subseteq Z_0 \subseteq Z$ (see Lemma 2.2), we have:
$$
(xy)z = x(yz) = x(zy) = (xz)y = y(xz) = (yx)z,
$$
whence $[x, y]z = (xy - yx)z = 0$. Further, since $xy, yx \in J^2 \subseteq Z_0 \subseteq Z$, the elements $xy$ and $yx$ are central in $R$. Therefore, by associativity of multiplication and centrality of $xy$, $yx$, we obtain:
$$
(xy)x = x(xy) = (xx)y = x^2 y = x(yx) = (yx)x,
$$
whence $[x, y]x = (xy - yx)x = 0$. Similarly $[x, y]y = 0$. Choose now elements $a, b \in J \setminus Z_0$ whose cosets generate $J/Z_0 \cong C_2 \times C_2$. Since $J = \langle a, b, Z_0 \rangle$ and $Z_0 \subseteq Z$, the commutator $[J, J]$ is generated by the single commutator $[a, b]$. By what was proved above, $[a, b]Z_0 = 0$, $[a, b]a = 0$, and $[a, b]b = 0$, whence $[a, b]J = 0$ and $[J, J] = \{0, [a, b]\} \subseteq W_0$. $\square$

\textbf{2.4. Lemma.}
In the ring $R$ the relations $JZ_0 \subseteq W_0$ and $Z_0^2 = 0$ hold.

\textbf{Proof.} By Remark 1.3 the order of the subgroup $W_0$ is $4$ or $2$. Consider both cases.

\textbf{Case 1:} $|W_0| = 4$. Choose an element $c \in Z_0 \setminus W_0$. Then $Z_0 = W_0 \cup (c + W_0)$. By Lemma 2.2, right multiplication by $c$ defines a homomorphism of additive groups $L_c : J \to Z_0$, $L_c(r) = rc$. If $r \in \ker L_c$, then $rc = 0$. Since $rW_0 = 0$ by the definition of $W_0$, the element $r$ annihilates the whole subgroup $Z_0 = \langle c, W_0 \rangle$, and by Lemma 2.2 we have $r \in Z_0$. Thus $W_0 \subseteq \ker L_c \subseteq Z_0$, whence $|\ker L_c| \in \{4, 8\}$. If $|\ker L_c| = 4$, then $|\mathrm{Im}(L_c)| = |J|/|\ker L_c| = 32/4 = 8 = |Z_0|$, that is, $L_c(J) = Z_0$, and there exists $u \in J$ such that $uc = c$. Then $(1 - u)c = 0$ and $c = 0$ --- a contradiction. Therefore $|\ker L_c| = 8$ and $\ker L_c = Z_0$. Since $Z_0 c = 0$ and $Z_0 W_0 = 0$, we obtain $Z_0^2 = 0$. Moreover, by Lemma 2.3, for any $x, y \in J$ we have $xy \in J^2 \subseteq Z_0 = \ker L_c$, whence $x(yc) = (xy)c = 0$. This means that $yc \in \mathrm{Ann}(J) \cap Z_0 = W_0$ for all $y \in J$, that is, $JZ_0 = Jc = W_0$.

\textbf{Case 2:} $|W_0| = 2$, $W_0 = \{0, w\}$. Denote $W_1 = \{z \in Z_0 \mid Jz \subseteq W_0\}$. Since the ideal $J$ is nilpotent and $JZ_0 \subseteq Z_0$, the induced action of $J$ on the nonzero quotient group $Z_0/W_0$ of order $4$ is nilpotent, so the subgroup of $J$-fixed elements $W_1/W_0$ is nontrivial, that is, $W_0 \subsetneq W_1 \subseteq Z_0$ and $|W_1| \in \{4, 8\}$.

Suppose first that $|W_1| = 4$, that is, $W_1 = \{0, w, c_1, c_1 + w\}$, and choose $c_2 \in Z_0 \setminus W_1$. By the nilpotency of the action of $J$ on $Z_0/W_1$ we have $Jc_2 \subseteq W_1$. The homomorphism $L_{c_1} : J \to W_0$ ($r \mapsto rc_1$) is nonzero, since $c_1 \notin W_0$, and has image $W_0$ of order $2$, so its kernel $H = \ker L_{c_1}$ is a subgroup of $J$ of order $32/2 = 16$. Note that every element $r \in H$ annihilates both $w \in W_0$ and $c_1 \in W_1$, and hence $HW_1 = 0$. For any $r \in H$ and $x \in J$, by Lemma 2.3 we have $[x, r] = xr - rx \in W_0$, whence $[x, r]c_2 = 0$ and
$$
x(rc_2) = (xr)c_2 = (rx)c_2 = r(xc_2) \in HW_1 = \{0\},
$$
since $xc_2 \in Jc_2 \subseteq W_1$. Therefore $J(rc_2) = 0$, that is, $rc_2 \in W_0$ for all $r \in H$. This means that the restriction $L_{c_2}|_H$ maps the subgroup $H$ of order $16$ into the group $W_0$ of order $2$. The kernel $K = \ker(L_{c_2}|_H) = \ker L_{c_1} \cap \ker L_{c_2}$ has order at least $16/2 = 8$ and annihilates all generators $w, c_1, c_2$ of the group $Z_0$. By Lemma 2.2 we have $K \subseteq \mathrm{Ann}_J(Z_0) \subseteq Z_0$. Comparing orders $|K| \ge 8 = |Z_0|$ gives $K = Z_0$, whence $Z_0 c_2 = 0$. But then for any $x, y \in J$, in view of $xy \in J^2 \subseteq Z_0$, we obtain $x(yc_2) = (xy)c_2 \in Z_0 c_2 = \{0\}$, whence $Jc_2 \subseteq W_0$ and $c_2 \in W_1$ --- a contradiction with the choice $c_2 \notin W_1$.

Thus the case $|W_1| = 4$ is impossible, and therefore $|W_1| = 8$, that is, $W_1 = Z_0$ and $JZ_0 \subseteq W_0$. Now consider the homomorphism of additive groups
$$
\Phi : J \to \mathrm{Hom}(Z_0/W_0, W_0), \qquad \Phi(r)(z + W_0) = rz.
$$
Since $|Z_0/W_0| = 4$ and $|W_0| = 2$, the order of the group $\mathrm{Hom}(Z_0/W_0, W_0)$ divides $4$. On the other hand, by Lemma 2.2 the kernel $\ker \Phi = \mathrm{Ann}_J(Z_0)$ is contained in $Z_0$. Hence
$$
8 = |Z_0| \ge |\ker \Phi| = \frac{|J|}{|\mathrm{Im}(\Phi)|} \ge \frac{32}{4} = 8,
$$
which implies $\ker \Phi = Z_0$, and hence $Z_0^2 = 0$. The lemma is completely proved. $\square$

\textbf{2.5. Completion of the proof of Theorem 1.4.}

\textbf{1.} The inclusions $J^2 \subseteq Z_0$, $[J, J] \subseteq W_0$, $JZ_0 \subseteq W_0$, and $Z_0^2 = 0$ are proved in Lemmas 2.3 and 2.4. From them it follows directly that $J^3 = J \cdot J^2 \subseteq JZ_0 \subseteq W_0$ and $J^4 = J \cdot J^3 \subseteq JW_0 = 0$. Since the ring $R$ is non-commutative, $[J, J] \neq 0$, whence $J^2 \neq 0$ and the nilpotency index of the radical $J$ equals $3$ or $4$.

\textbf{2.} Since $R/J \cong \mathbb{F}_2$ and $1 \in Z$, the element $2\cdot 1 = 1 + 1$ belongs to $Z \cap J = Z_0$. Since $Z_0^2 = 0$, for any element $z \in Z_0$ we have $2z = (2\cdot 1)z \in Z_0^2 = \{0\}$. Hence $2Z_0 = 0$ and $Z_0 \cong \mathbb{F}_2^3$. In particular, $4\cdot 1 = 2(2\cdot 1) \in 2Z_0 = \{0\}$, whence $\mathrm{char}(R) = 2$ (when $2\cdot 1 = 0$) or $\mathrm{char}(R) = 4$ (when $2\cdot 1 \neq 0$). Since $J^2 \subseteq Z_0$, from $2Z_0 = 0$ it follows immediately that $2J^2 = 0$.

\textbf{3.} Let $\dim_{\mathbb{F}_2} W_0 = 2$, $W_0 = \langle s_1, s_2 \rangle_{\mathbb{F}_2}$, and $Z_0 = \langle c, s_1, s_2 \rangle_{\mathbb{F}_2}$. As shown in the proof of Lemma 2.4 (Case 1), the map $L_c : J \to W_0$ ($r \mapsto rc$) has kernel $Z_0$ and induces an isomorphism $\overline{L}_c : J/Z_0 \xrightarrow{\sim} W_0$. Choose elements $a, b \in J \setminus Z_0$ such that $ac = s_1$ and $bc = s_2$. Since $J^2 \subseteq Z_0$, for any $x, y \in \{a, b\}$ we have the decomposition
$$
xy \equiv \gamma(x, y)c \pmod{W_0}, \qquad \gamma(x, y) \in \mathbb{F}_2,
$$
where $\gamma(a, b) = \gamma(b, a)$ in view of $[a, b] \in W_0$. For any $y \in \{a, b\}$, by the associativity law $(ay)b = a(yb)$. Taking into account $W_0 J = JW_0 = 0$ and the centrality of $c$, we obtain:
$$
(ay)b = \gamma(a, y)cb = \gamma(a, y)s_2, \qquad a(yb) = \gamma(y, b)ac = \gamma(y, b)s_1.
$$
In view of the linear independence of the vectors $s_1, s_2 \in W_0$ over $\mathbb{F}_2$, the equality $\gamma(a, y)s_2 = \gamma(y, b)s_1$ implies $\gamma(a, y) = \gamma(y, b) = 0$ for all $y \in \{a, b\}$. Consequently, $\gamma \equiv 0$, that is, $J^2 \subseteq W_0$ and $J^3 = 0$.

\section{Classification in characteristic 2 ($\mathbb{F}_2$-algebras)}

Let $\mathrm{char}(R) = 2$, that is, $2\cdot 1 = 0$ and $R$ is a six-dimensional algebra over the field $\mathbb{F}_2$. By Theorem 1.4 there are exactly three mutually exclusive cases:
\begin{itemize}
    \item \textbf{Family 1:} $\dim_{\mathbb{F}_2} W_0 = 2$, with $J^2 \subseteq W_0$ and $J^3 = 0$;
    \item \textbf{Family 2:} $\dim_{\mathbb{F}_2} W_0 = 1$ and $J^2 \subseteq W_0$, with $J^3 = 0$;
    \item \textbf{Family 3:} $\dim_{\mathbb{F}_2} W_0 = 1$ and $J^2 \nsubseteq W_0$, with $J^3 = W_0 \neq 0$, $J^4 = 0$.
\end{itemize}
We investigate each of these families and find all isomorphism classes.

\textbf{Family 1 ($\dim_{\mathbb{F}_2} W_0 = 2$)}

\textbf{3.1. Lemma.}
Algebras of Family 1 have a basis $\{1, a, b, c, s_1, s_2\}$, $J(R) = \langle a, b, c, s_1, s_2 \rangle_{\mathbb{F}_2}$, $Z(R) = \langle 1, c, s_1, s_2 \rangle_{\mathbb{F}_2}$, and $W_0 = \langle s_1, s_2 \rangle_{\mathbb{F}_2}$. The multiplication of basis elements of the radical is given by the rules:
$$
ac = s_1; \quad bc = s_2; \quad ab = s_1, \,\, ba = 0; \quad c^2 = 0. \eqno (3.1)
$$
The elements $s_1, s_2$ annihilate $J(R)$. The squares $a^2 = \alpha$ and $b^2 = \beta$ are chosen from $W_0$, which parametrizes the algebras within the family.

\textbf{Proof.} Let $W_0 = \langle s_1, s_2 \rangle$ and $Z_0 = \langle c, s_1, s_2 \rangle$. Consider the right multiplication operator $L_c: J \to W_0$ ($r \mapsto rc$). Since $\dim J = 5$ and $\dim W_0 = 2$, the dimension of the kernel $\dim \ker L_c \ge 3$. Moreover, $rZ_0\neq 0$ for a noncentral element $r$, otherwise the ring $R$ would not be centrally essential. Then $\ker L_c \subseteq Z_0$. Taking into account $\dim Z_0 = 3$, we obtain $\ker L_c = Z_0$, which means $Z_0 c = 0$. Hence $Z^2_0 = 0$ in view of $Z_0W_0 = 0$ (in particular, $c^2 = 0$). Consequently, $L_c$ induces a well-defined isomorphism $J/Z_0 \to W_0$. Choose elements $a, b \in J \setminus Z_0$ whose cosets form a basis of $J/Z_0$, such that $ac = s_1$ and $bc = s_2$.

Since the algebra $R$ is non-commutative, the commutator $[a,b] \neq 0$. Since $J^2 \subseteq W_0$, we have $t = [a,b] \in W_0 \setminus \{0\}$. Let us show that by a suitable change of basis in $J$ one can achieve the equality $[a,b] = s_1$. If $t = s_1$, the required equality holds. Let us show that for $t = s_2$ one can achieve the condition $[a', b'] = s_1'$. Perform the change $a' = b$, $b' = a + b$. Since the characteristic of the ring $R$ equals $2$, we have:
$$
[a', b'] = [b, a + b] = [b, a] = -[a, b] = s_2.
$$
Next, pass to a new basis of the space $W_0$, setting $s_1' = s_2$ and $s_2' = s_1 + s_2$. In the new basis we obtain the required relation $[a', b'] = s_1'$. Multiplication by $c$ is also consistent with the new basis: $a' c = b c = s_2 = s_1'$ and $b' c = (a + b) c = s_1 + s_2 = s_2'$. In the case $t = s_1 + s_2$, by the changes $a' = a + b$, $b' = a$ and $s_1' = s_1 + s_2$, $s_2' = s_1$, we again obtain $[a', b'] = s_1'$. Thus we may assume that $[a,b] = s_1$.

Also, without loss of generality, we may set $ba = 0$. Indeed, taking into account $Z^2_0 = 0$, under the change $a' = a + x, b' = b + y$ (where $x,y \in Z_0$) we have $b'a' = ba + xb + ya$. Since the operators $L_a, L_b$ map $Z_0$ onto the one-dimensional subspaces $\langle s_1 \rangle$ and $\langle s_2 \rangle$, respectively, the elements of $W_0$ are linear combinations of $xb + ya$. Hence there exist $x,y \in Z_0$ such that $xb + ya = ba$ and $b'a' = 0$. Moreover, $a'c = (a+x)c = ac + xc = s_1 + 0 = s_1$, since $xc = 0$ for $x \in Z_0$. Further, $b'c = (b + y)c = bc + yc = s_2 + 0 = s_2$. Consequently, the change is correct. The squares $a^2 = \alpha$ and $b^2 = \beta$ are chosen from $W_0$. Different pairs $(\alpha, \beta)$ parametrize this family.

Note that under such a change the commutator $[a', b'] = [a, b] = s_1$ is preserved. Finally, since the squares $a^2 = \alpha$ and $b^2 = \beta$ belong to $W_0$, they do not affect the commutators of basis elements or the multiplication by the element $c$. Therefore, as in the proof of Lemma 2.3, for any $\alpha, \beta \in W_0$ the equalities $Z(R) = \langle 1, c, s_1, s_2 \rangle_{\mathbb{F}_2}$, $W_0 = \langle s_1, s_2 \rangle_{\mathbb{F}_2}$ hold, and for any noncentral element $r \in R \setminus Z(R)$ we have $0 \neq rc \in Z(R)$, i.e., the algebra $R$ is centrally essential. $\square$

Denote by $R_1(\alpha, \beta)$ the algebra with parameters $(\alpha, \beta) \in W_0 \times W_0$. In total there are $4 \times 4 = 16$ pairs of parameters $(\alpha, \beta) \in W_0^2$.

\textbf{3.2. Proposition.}
The sixteen algebras $R_1(\alpha, \beta)$ ($\alpha, \beta \in W_0$) of Family 1 split into exactly $10$ pairwise non-isomorphic classes. The complete list of isomorphism classes and the corresponding orbits of pairs $(\alpha, \beta)$ is given in Table 1.

\textbf{Proof.} Let $R = R_1(\alpha, \beta)$ be the algebra with basis $\{1, a, b, c, s_1, s_2\}$ satisfying relations (3.1). Let us describe all other bases $\{1, a', b', c', s_1', s_2'\}$ of the algebra $R$ in which the multiplication table has the same canonical form (3.1) with some parameters $\alpha' = (a')^2, \beta' = (b')^2 \in W_0$.

First, the subspaces $J$, $Z_0 = Z(R) \cap J$, and $W_0 = \mathrm{Ann}(J) \cap Z_0$ are invariant under any automorphism (and isomorphism) of algebras. Moreover, the commutator $[J, J] = \{0, [a, b]\} = \{0, s_1\}$ uniquely singles out the nonzero vector $s_1 \in W_0$, so in any canonical basis the equality $s_1' = [a', b'] = s_1$ must hold.

Second, any element $c' \in Z_0 \setminus W_0$ has the form $c' = c + w_c$, where $w_c \in W_0$. Since $JW_0 = 0$, for any $r \in J$ we have $rc' = r(c + w_c) = rc$. Hence the linear isomorphism $\overline{L}_c : J/Z_0 \xrightarrow{\sim} W_0$, $r + Z_0 \mapsto rc$, does not depend on the choice of the element $c' \in Z_0 \setminus W_0$ and is an internal invariant of the algebra $R$. In particular, the condition $a'c' = s_1' = s_1$ uniquely determines the coset of the element $a'$ modulo $Z_0$: $a' \equiv a \pmod{Z_0}$.

Third, since $\{a' + Z_0, b' + Z_0\}$ is a basis of the two-dimensional $\mathbb{F}_2$-space $J/Z_0$, for the coset $b' + Z_0$ there are exactly two possibilities:
$$
b' \equiv b \pmod{Z_0} \quad \text{or} \quad b' \equiv a + b \pmod{Z_0}.
$$
Fourth, note that the squaring map $\mathrm{Sq} : J \to W_0$, $r \mapsto r^2$, is constant on cosets modulo the subspace $Z_0$. Indeed, for any $z \in Z_0$, in view of $z \in Z(R)$, $\mathrm{char}(R) = 2$, and $Z_0^2 = 0$, we have $(r + z)^2 = r^2 + 2rz + z^2 = r^2$. Consequently, the elements $(a')^2$ and $(b')^2$ depend only on the choice of the cosets $a' + Z_0$ and $b' + Z_0$.
\begin{itemize}
    \item If $b' \equiv b \pmod{Z_0}$, then $s_1' = a'c' = ac = s_1$, $s_2' = b'c' = bc = s_2$, $(a')^2 = a^2 = \alpha$, and $(b')^2 = b^2 = \beta$, that is, the parameters $(\alpha, \beta)$ do not change.
    \item If $b' \equiv a + b \pmod{Z_0}$, set $s_1' = s_1$, $s_2' = (a + b)c = s_1 + s_2$, $c' = c$, $a' = a + x$, and $b' = a + b + y$, where $x, y \in Z_0$ are chosen so that
$$
b'a' = (a + b + y)(a + x) = a^2 + ba + ya + x(a + b) = \alpha + ya + x(a + b) = 0.
$$
This is always possible and unique up to $W_0$, since $ca = s_1$ and $c(a + b) = s_1 + s_2$ form a basis of $W_0$. Then in the new canonical basis $\{1, a', b', c', s_1', s_2'\}$ all relations (3.1) hold, and the squares of the elements $a'$ and $b'$ equal:
$$
\alpha' = (a')^2 = a^2 = \alpha, \quad \beta' = (b')^2 = (a + b)^2 = a^2 + b^2 + ab + ba = \alpha + \beta + s_1.
$$
\end{itemize}
Write the coordinates of the vectors $\alpha, \beta \in W_0$ in the original basis $\{s_1, s_2\}$:
$$
\alpha = \alpha_1 s_1 + \alpha_2 s_2, \quad \beta = \beta_1 s_1 + \beta_2 s_2 \quad (\alpha_1, \alpha_2, \beta_1, \beta_2 \in \mathbb{F}_2).
$$
Passing to the new basis $\{s_1' = s_1, s_2' = s_1 + s_2\}$ of the subspace $W_0$ gives $s_2 = s_1' + s_2'$. Therefore the coordinates of the vectors $\alpha' = \alpha$ and $\beta' = \alpha + \beta + s_1$ in the basis $\{s_1', s_2'\}$ equal:
$$
\alpha' = \alpha_1 s_1' + \alpha_2(s_1' + s_2') = (\alpha_1 + \alpha_2)s_1' + \alpha_2 s_2',
$$
$$
\beta' = (\alpha_1 + \beta_1 + 1)s_1' + (\alpha_2 + \beta_2)(s_1' + s_2') =
$$
$$
(\alpha_1 + \alpha_2 + \beta_1 + \beta_2 + 1)s_1' + (\alpha_2 + \beta_2)s_2'.
$$
Thus two algebras $R_1(\alpha, \beta)$ and $R_1(\alpha', \beta')$ are isomorphic if and only if their coordinate quadruples $(\alpha_1, \alpha_2, \beta_1, \beta_2) \in \mathbb{F}_2^4$ lie in the same orbit under the involution $\tau : \mathbb{F}_2^4 \to \mathbb{F}_2^4$:
$$
\tau(\alpha_1, \alpha_2, \beta_1, \beta_2) = (\alpha_1 + \alpha_2, \, \alpha_2, \, \alpha_1 + \alpha_2 + \beta_1 + \beta_2 + 1, \, \alpha_2 + \beta_2).
$$
Let us find the fixed points of the involution $\tau$. The system of equations $\tau(\alpha_1, \alpha_2, \beta_1, \beta_2) = (\alpha_1, \alpha_2, \beta_1, \beta_2)$ is equivalent to the conditions
$$
\alpha_1 + \alpha_2 = \alpha_1, \qquad \alpha_1 + \alpha_2 + \beta_1 + \beta_2 + 1 = \beta_1, \qquad \alpha_2 + \beta_2 = \beta_2,
$$
whence $\alpha_2 = 0$ and $\beta_2 = \alpha_1 + 1$, while $\alpha_1, \beta_1 \in \mathbb{F}_2$ are arbitrary. This gives exactly $4$ fixed points (orbits of length $1$), and the remaining $16 - 4 = 12$ quadruples split into $12/2 = 6$ orbits of length $2$. In total we obtain $4 + 6 = 10$ isomorphism classes, listed in Table 1. $\square$

\begin{table}[htbp]
\centering
\caption{Classification of algebras of Family 1 ($\operatorname{char}(R) = 2$, $\dim_{\mathbb{F}_2} W_0 = 2$, $J^3 = 0$)}
\vspace{0.5em}
\begin{tabular}{ |c|c|c|c|}
\hline
Class & Representative $(\alpha, \beta)$ & Orbit of pairs $(a^2, b^2)$ in the basis $\{s_1, s_2\}$ & Orbit length \\
\hline
$R_{1.1}$ & $(0, s_2)$ & $\{(0, s_2)\}$ & $1$ \\
\hline
$R_{1.2}$ & $(0, s_1 + s_2)$ & $\{(0, s_1 + s_2)\}$ & $1$ \\
\hline
$R_{1.3}$ & $(s_1, 0)$ & $\{(s_1, 0)\}$ & $1$ \\
\hline
$R_{1.4}$ & $(s_1, s_1)$ & $\{(s_1, s_1)\}$ & $1$ \\
\hline
$R_{1.5}$ & $(0, 0)$ & $\{(0, 0), (0, s_1)\}$ & $2$ \\
\hline
$R_{1.6}$ & $(s_1, s_2)$ & $\{(s_1, s_2), (s_1, s_1 + s_2)\}$ & $2$ \\
\hline
$R_{1.7}$ & $(s_2, 0)$ & $\{(s_2, 0), (s_1 + s_2, s_2)\}$ & $2$ \\
\hline
$R_{1.8}$ & $(s_2, s_1)$ & $\{(s_2, s_1), (s_1 + s_2, s_1 + s_2)\}$ & $2$ \\
\hline
$R_{1.9}$ & $(s_2, s_2)$ & $\{(s_2, s_2), (s_1 + s_2, s_1)\}$ & $2$ \\
\hline
$R_{1.10}$ & $(s_2, s_1 + s_2)$ & $\{(s_2, s_1 + s_2), (s_1 + s_2, 0)\}$ & $2$ \\
\hline
\end{tabular}
\end{table}

\textbf{Family 2 ($\dim_{\mathbb{F}_2} W_0 = 1$ and $J^2 \subseteq W_0$)}

\textbf{3.3. Lemma.} An algebra of Family 2 ($\dim_{\mathbb{F}_2} W_0 = 1$) has a basis $\{1, a, b, c_1, c_2, w\}$, $J(R) = \langle a, b, c_1, c_2, w \rangle_{\mathbb{F}_2}$, $Z(R) = \langle 1, c_1, c_2, w \rangle_{\mathbb{F}_2}$, and $W_0 = \langle w \rangle_{\mathbb{F}_2}$. The multiplication of basis elements of the radical is given by the rules:
$$
ac_1 = w, ac_2 = 0; bc_1 = 0, bc_2 = w; ab = w, ba = 0; c_1^2 = c_2^2 = c_1 c_2 = 0. \eqno (3.2)
$$
The element $w$ annihilates the radical $J(R)$. The squares $a^2 = \alpha w$ and $b^2 = \beta w$ are parametrized by scalars $\alpha, \beta \in \mathbb{F}_2$.

\textbf{Proof.} Let $W_0 = \langle w \rangle$ and $Z_0 = \langle c_1, c_2, w \rangle$. Since $J^2 \subseteq W_0$ and $JW_0 = 0$, the linear map
$$
\Phi: J \to \mathrm{Hom}_{\mathbb{F}_2}(Z_0/W_0, W_0), \qquad \Phi(r)(z + W_0) = rz
$$
is well defined. If $r \in \ker \Phi$, then $rZ_0 = 0$. From the central essentiality of the algebra $R$ it follows that no noncentral element $r \in J \setminus Z_0$ can annihilate $Z_0$, whence $\ker \Phi \subseteq Z_0$. Since $\dim J = 5$ and $\dim \mathrm{Hom}_{\mathbb{F}_2}(Z_0/W_0, W_0) = 2$, by the rank-nullity theorem we have $\dim \ker \Phi \ge 5 - 2 = 3$. Taking into account $\dim Z_0 = 3$, we obtain $\ker \Phi = Z_0$, which means $Z_0^2 = 0$, i.e., $c_1^2 = c_2^2 = c_1c_2 = 0$. Moreover, $\Phi$ induces an isomorphism of two-dimensional spaces $\bar{\Phi}: J/Z_0 \to \mathrm{Hom}_{\mathbb{F}_2}(Z_0/W_0, W_0)$. Consequently, one can choose elements $a, b \in J \setminus Z_0$ whose images in $J/Z_0$ form a basis dual to the basis $\{c_1 + W_0, c_2 + W_0\}$ of the space $Z_0/W_0$, so that
$$
ac_1 = w, \quad ac_2 = 0, \quad bc_1 = 0, \quad bc_2 = w.
$$
Since the algebra $R$ is non-commutative, the commutator $[a,b] \neq 0$. Since $[a,b] \in J^2 \subseteq W_0 = \langle w \rangle$, we have $ab + ba = [a,b] = w$. Since $ab, ba \in \{0, w\}$, two cases are possible: $ab = w, ba = 0$ or $ab = 0, ba = w$. In the second case perform the change $a' = a + c_2$. Since $Z_0^2 = 0$, the relations $a'c_1 = w$ and $a'c_2 = 0$ are preserved, while $a'b = ab + c_2b = w$ and $ba' = ba + bc_2 = 0$. Thus, without loss of generality, we may assume that $ab = w$ and $ba = 0$. The squares $a^2 = \alpha w$ and $b^2 = \beta w$ are parametrized by scalars $\alpha, \beta \in \mathbb{F}_2$.

Conversely, for any pair $(\alpha, \beta) \in \mathbb{F}_2^2$, the algebra $R$ defined by these rules has center $Z(R) = \langle 1, c_1, c_2, w \rangle_{\mathbb{F}_2}$ and $W_0 = \langle w \rangle_{\mathbb{F}_2}$. For an arbitrary noncentral element $r = k_0 1 + k_a a + k_b b + z \in R \setminus Z(R)$ ($z \in Z_0$, $(k_a, k_b) \neq (0,0)$) we have $rc_1 = k_0 c_1 + k_a w$ and $rc_2 = k_0 c_2 + k_b w$. Since $(k_a, k_b) \neq (0,0)$, at least one of the products $rc_1, rc_2$ is nonzero and belongs to $Z_0 \subseteq Z(R)$, whence the central essentiality of $R$ follows. $\square$

\textbf{3.4. Proposition.} Family 2 consists of exactly $2$ pairwise non-isomorphic algebras, corresponding to the parameter pairs $(\alpha, \beta) = (0, 0)$ (class $R_{2.1}$) and $(\alpha, \beta) = (1, 1)$ (class $R_{2.2}$).

\textbf{Proof.} For an arbitrary element $x = x_1 a + x_2 b + z \in J$ ($x_1, x_2 \in \mathbb{F}_2$, $z \in Z_0$), in view of $Z_0^2 = 0$ and $\mathrm{char}(R) = 2$, we have:
$$
x^2 = x_1^2 a^2 + x_2^2 b^2 + x_1 x_2 (ab + ba) = (x_1 \alpha + x_2 \beta + x_1 x_2)w = Q(x_1, x_2)w.
$$
The quadratic form $Q(x_1, x_2) = \alpha x_1^2 + x_1 x_2 + \beta x_2^2$ over the field $\mathbb{F}_2$ has invariant $\mathrm{Arf}(Q) = \alpha\beta \in \mathbb{F}_2$:
\begin{itemize}
    \item For $(\alpha, \beta) \in \{(0, 0), (1, 0), (0, 1)\}$ ($\mathrm{Arf}(Q) = 0$) the form vanishes on exactly $3$ pairs $(x_1, x_2) \in \mathbb{F}_2^2$, so the radical $J$ contains exactly $3 \cdot |Z_0| = 24$ elements with zero square. Moreover, the algebra with $(\alpha, \beta) = (1, 0)$ is transformed into the algebra with $(0, 0)$ by the change $a' = b + c_1$, $b' = a + b$, $c_1' = c_1 + c_2$, $c_2' = c_1$, and the algebra with $(0, 1)$ is transformed into $(0, 0)$ by the change $a' = a$, $b' = a + b + c_2$, $c_1' = c_1 + c_2$, $c_2' = c_2$.
    \item For $(\alpha, \beta) = (1, 1)$ ($\mathrm{Arf}(Q) = 1$) the form $Q(x_1, x_2) = x_1^2 + x_1 x_2 + x_2^2$ is anisotropic over $\mathbb{F}_2$ and vanishes only for $x_1 = x_2 = 0$, so the radical $J$ contains exactly $1 \cdot |Z_0| = 8$ elements with zero square.
\end{itemize}
Consequently, Family 2 contains exactly $2$ isomorphism classes. $\square$

\textbf{Family 3 ($\dim_{\mathbb{F}_2} W_0 = 1$ and $J^2 \nsubseteq W_0$)}

We now turn to the description of the class of $\mathbb{F}_2$-algebras not satisfying the condition $J^2 \subseteq W_0$. By Theorem 1.4, in this case $\dim_{\mathbb{F}_2} W_0 = 1$, $W_0 = \langle w \rangle_{\mathbb{F}_2}$, $Z_0 = \langle c_1, c_2, w \rangle_{\mathbb{F}_2}$, $JZ_0 \subseteq W_0$, and $Z_0^2 = 0$. Since the kernel of the map $\Phi : J \to \mathrm{Hom}_{\mathbb{F}_2}(Z_0/W_0, W_0)$ coincides with $Z_0$ by Lemma 2.4, we can still choose a basis $\{c_1 + W_0, c_2 + W_0\}$ in $Z_0/W_0$ and a dual basis $\{a + Z_0, b + Z_0\}$ in $J/Z_0$ so that:
$$
ac_1 = w, \,\, ac_2 = 0, \,\, bc_1 = 0, \,\, bc_2 = w, \,\, Z_0^2 = 0, \,\, JW_0 = 0. \eqno (3.3)
$$
\textbf{3.5. Lemma.}
Let the subspace $J = \langle a, b, c_1, c_2, w \rangle_{\mathbb{F}_2}$ satisfy conditions (3.3), with $J^2 \subseteq Z_0$ and $[a, b] = w$. Multiplication on $J$ is associative if and only if, for a suitable choice of the representative $a$ in the coset $a + Z_0$, the products of the elements $a$ and $b$ have the form:
$$
a^2 = \lambda_1 c_1 + \lambda_2 c_2 + \alpha w, \quad ab = \lambda_2 c_1 + \lambda_3 c_2 + w,
$$
$$
ba = \lambda_2 c_1 + \lambda_3 c_2, \quad b^2 = \lambda_3 c_1 + \lambda_4 c_2 + \beta w, \eqno (3.4)
$$
where $\lambda_1, \lambda_2, \lambda_3, \lambda_4, \alpha, \beta \in \mathbb{F}_2$. Moreover, the condition $J^2 \nsubseteq W_0$ (that is, $J^3 = W_0 \neq 0$) is equivalent to $(\lambda_1, \lambda_2, \lambda_3, \lambda_4) \neq (0, 0, 0, 0)$.

\textbf{Proof.} Since $J^2 \subseteq Z_0 = \langle c_1, c_2, w \rangle_{\mathbb{F}_2}$ and $ab + ba = [a, b] = w$, the products of the elements $a, b$ in general form are written as
$$
a^2 = \gamma_{11}^{(1)} c_1 + \gamma_{11}^{(2)} c_2 + \alpha w, \quad ab = \gamma_{12}^{(1)} c_1 + \gamma_{12}^{(2)} c_2 + \delta w,
$$
$$
ba = \gamma_{12}^{(1)} c_1 + \gamma_{12}^{(2)} c_2 + (\delta + 1)w, \quad b^2 = \gamma_{22}^{(1)} c_1 + \gamma_{22}^{(2)} c_2 + \beta w.
$$
If $\delta = 0$, then the change $a \mapsto a + c_2$ preserves relations (3.3), does not change $a^2 = (a + c_2)^2$, and adds $c_2 b = w$ to $ab$ and $ba$, making the coefficient of $w$ in $ba$ zero ($\delta = 1$).

Let us verify the associativity condition $(xy)z = x(yz)$ on basis elements $x, y, z \in \{a, b, c_1, c_2, w\}$. If at least one of the elements $x, y, z$ belongs to $Z_0 = \langle c_1, c_2, w \rangle_{\mathbb{F}_2}$, then both sides of the equality $(xy)z = x(yz)$ vanish, since $J^2 Z_0 \subseteq Z_0^2 = 0$, $J(JZ_0) \subseteq JW_0 = 0$, and $(JZ_0)J \subseteq W_0 J = 0$. Let now $x, y, z \in \{a, b\}$. Since $xy \in Z_0 \subseteq Z(R)$ and $xy - yx \in W_0 \subseteq \mathrm{Ann}(J)$, the identities $(a^2)a = a(a^2)$, $(b^2)b = b(b^2)$, $(ab)a = a(ba)$, and $(ba)b = b(ab)$ hold automatically. Four triples remain. For $(a, a, b)$ and $(b, a, a)$ we have:
$$
(a^2)b = (\gamma_{11}^{(1)} c_1 + \gamma_{11}^{(2)} c_2 + \alpha w)b = \gamma_{11}^{(2)} w, \,\, a(ab) = a(\gamma_{12}^{(1)} c_1 + \gamma_{12}^{(2)} c_2 + w) = \gamma_{12}^{(1)} w,
$$
whence $(a^2)b = a(ab) \iff \gamma_{11}^{(2)} = \gamma_{12}^{(1)}$.
For $(a, b, b)$ and $(b, b, a)$:
$$
(ab)b = (\gamma_{12}^{(1)} c_1 + \gamma_{12}^{(2)} c_2 + w)b = \gamma_{12}^{(2)} w, \,\, a(b^2) = a(\gamma_{22}^{(1)} c_1 + \gamma_{22}^{(2)} c_2 + \beta w) = \gamma_{22}^{(1)} w,
$$
whence $(ab)b = a(b^2) \iff \gamma_{12}^{(2)} = \gamma_{22}^{(1)}$.

Setting $\lambda_1 = \gamma_{11}^{(1)}$, $\lambda_2 = \gamma_{11}^{(2)} = \gamma_{12}^{(1)}$, $\lambda_3 = \gamma_{12}^{(2)} = \gamma_{22}^{(1)}$, and $\lambda_4 = \gamma_{22}^{(2)}$, we obtain relations (3.4). Moreover, for every $r = k_a a + k_b b + z \in J \setminus Z_0$ ($(k_a, k_b) \neq (0, 0)$) we have $rc_1 = k_a w$ and $rc_2 = k_b w$, so that $\mathrm{Ann}_J(Z_0) = Z_0$ and, by Lemma 2.2, the algebra $R$ is centrally essential. $\square$

\textbf{3.6. Proposition.}
Family 3 ($\mathrm{char}(R) = 2$, $\dim_{\mathbb{F}_2} W_0 = 1$, $J^3 \neq 0$) consists of exactly $7$ pairwise non-isomorphic algebras, listed in Table 2.

\textbf{Proof.} On the quotient space $V = J/Z_0$ with basis $\{\overline{a} = a + Z_0, \overline{b} = b + Z_0\}$, the multiplication in the radical $J$ defines an internally invariant symmetric trilinear form with values in $W_0 = \langle w \rangle_{\mathbb{F}_2}$:
$$
T : V \times V \times V \to W_0, \qquad T(u + Z_0, v + Z_0, t + Z_0) = (uv)t.
$$
In the basis $\{\overline{a}, \overline{b}\}$ the values of the form $T$ equal:
$$
T(\overline{a}, \overline{a}, \overline{a}) = \lambda_1 w, \quad T(\overline{a}, \overline{a}, \overline{b}) = \lambda_2 w, \quad T(\overline{a}, \overline{b}, \overline{b}) = \lambda_3 w, \quad T(\overline{b}, \overline{b}, \overline{b}) = \lambda_4 w.
$$
Under a change of basis of the space $V = J/Z_0$ by a nonsingular matrix $M = \begin{pmatrix} p & q \\ r & s \end{pmatrix} \in \mathrm{GL}_2(\mathbb{F}_2)$ ($\overline{a}' = p\overline{a} + q\overline{b}$, $\overline{b}' = r\overline{a} + s\overline{b}$) and the corresponding dual change of basis in $Z_0/W_0$, the new structure constants $(\lambda_1', \lambda_2', \lambda_3', \lambda_4')$ are computed by the formulas:
$$
\lambda_1' w = T(\overline{a}', \overline{a}', \overline{a}') = (\lambda_1 p^3 + \lambda_2 p^2 q + \lambda_3 p q^2 + \lambda_4 q^3)w,
$$
$$
\lambda_2' w = T(\overline{a}', \overline{a}', \overline{b}') = (\lambda_1 p^2 r + \lambda_2 p^2 s + \lambda_3 q^2 r + \lambda_4 q^2 s)w,
$$
$$
\lambda_3' w = T(\overline{a}', \overline{b}', \overline{b}') = (\lambda_1 p r^2 + \lambda_2 q r^2 + \lambda_3 p s^2 + \lambda_4 q s^2)w,
$$
$$
\lambda_4' w = T(\overline{b}', \overline{b}', \overline{b}') = (\lambda_1 r^3 + \lambda_2 r^2 s + \lambda_3 r s^2 + \lambda_4 s^3)w.
$$
A direct expansion of the brackets in the polynomial ring $\mathbb{F}_2[X, Y]$, taking into account $(pX + rY)^2 = p^2 X^2 + r^2 Y^2$, shows that the quadruple $(\lambda_1', \lambda_2', \lambda_3', \lambda_4')$ coincides exactly with the coefficient vector of the transformed homogeneous binary cubic form $P'(X, Y) = P(pX + rY, qX + sY)$, where $P(X, Y) = \lambda_1 X^3 + \lambda_2 X^2 Y + \lambda_3 X Y^2 + \lambda_4 Y^3$.

The fifteen nonzero cubic forms in $\mathbb{F}_2[X, Y]_3 \setminus \{0\}$ under the action of the group $\mathrm{GL}_2(\mathbb{F}_2)$ of order $6$ split into exactly $5$ orbits, uniquely distinguished by the number of roots in the projective line $\mathbb{P}^1(\mathbb{F}_2)$ and by decomposability into linear factors:
\begin{enumerate}
    \item \textbf{Orbit 1} ($3$ forms, triple root in $\mathbb{P}^1(\mathbb{F}_2)$): representative $P(X, Y) = X^3$, $(\lambda_1, \lambda_2, \lambda_3, \lambda_4) = (1, 0, 0, 0)$.
    \item \textbf{Orbit 2} ($6$ forms, double and simple roots in $\mathbb{P}^1(\mathbb{F}_2)$): representative $P(X, Y) = X^2 Y$, $(\lambda_1, \lambda_2, \lambda_3, \lambda_4) = (0, 1, 0, 0)$.
    \item \textbf{Orbit 3} ($1$ form, three distinct roots in $\mathbb{P}^1(\mathbb{F}_2)$): representative $P(X, Y) = X^2 Y + X Y^2$, $(\lambda_1, \lambda_2, \lambda_3, \lambda_4) = (0, 1, 1, 0)$.
    \item \textbf{Orbit 4} ($3$ forms, product of a linear and an irreducible quadratic form): representative $P(X, Y) = X(X^2 + XY + Y^2) = X^3 + X^2 Y + X Y^2$, $(\lambda_1, \lambda_2, \lambda_3, \lambda_4) = (1, 1, 1, 0)$.
    \item \textbf{Orbit 5} ($2$ irreducible cubic forms): representative $P(X, Y) = X^3 + X^2 Y + Y^3$, $(\lambda_1, \lambda_2, \lambda_3, \lambda_4) = (1, 1, 0, 1)$.
\end{enumerate}

Let us now investigate how the parameters $(\alpha, \beta) \in \mathbb{F}_2^2$ transform for a fixed representative $(\lambda_1, \lambda_2, \lambda_3, \lambda_4)$. For any choice of representatives $c_1' = c_1 + \mu_1 w$, $c_2' = c_2 + \mu_2 w$ ($\mu_1, \mu_2 \in \mathbb{F}_2$) of the cosets $c_1 + W_0, c_2 + W_0$ and for the shift $a' = a + (\lambda_2 \mu_1 + \lambda_3 \mu_2)c_2'$ preserving the equality $b'a' = \lambda_2 c_1' + \lambda_3 c_2'$, in view of $(r + z)^2 = r^2$ ($z \in Z_0$) we obtain:
$$
(a')^2 = a^2 = \lambda_1 c_1' + \lambda_2 c_2' + (\alpha + \lambda_1 \mu_1 + \lambda_2 \mu_2)w,
$$
$$
(b')^2 = b^2 = \lambda_3 c_1' + \lambda_4 c_2' + (\beta + \lambda_3 \mu_1 + \lambda_4 \mu_2)w.
$$
Consequently, the change $(c_1, c_2) \mapsto (c_1 + \mu_1 w, c_2 + \mu_2 w)$ shifts the vector $(\alpha, \beta)^T \in \mathbb{F}_2^2$ by the vector $\Lambda(\mu_1, \mu_2)^T$, where $\Lambda = \begin{pmatrix} \lambda_1 & \lambda_2 \\ \lambda_3 & \lambda_4 \end{pmatrix}$.
\begin{itemize}
    \item For Orbits 3, 4, and 5 the matrices $\Lambda$ have the form $\begin{pmatrix} 0 & 1 \\ 1 & 0 \end{pmatrix}$, $\begin{pmatrix} 1 & 1 \\ 1 & 0 \end{pmatrix}$, and $\begin{pmatrix} 1 & 1 \\ 0 & 1 \end{pmatrix}$, respectively. All of them are nonsingular ($\det \Lambda = 1 \in \mathbb{F}_2$), so by the choice of $(\mu_1, \mu_2) \in \mathbb{F}_2^2$ any pair $(\alpha, \beta)$ is reduced to $(0, 0)$. Thus Orbits 3, 4, and 5 give one isomorphism class each ($R_{3.5}, R_{3.6}, R_{3.7}$).
    \item For Orbit 1 ($\Lambda = \begin{pmatrix} 1 & 0 \\ 0 & 0 \end{pmatrix}$) and Orbit 2 ($\Lambda = \begin{pmatrix} 0 & 1 \\ 0 & 0 \end{pmatrix}$) the rank of the matrix $\Lambda$ equals $1$. Therefore, by the shift $(\mu_1, \mu_2)$ the parameter $\alpha$ is reduced to $0$, while the parameter $\beta \in \{0, 1\}$ is invariant. Indeed, in both cases $b + Z_0$ is the unique nonzero coset in $J/Z_0$ whose square lies in $W_0$ ($b^2 = \beta w$), and $(b + z)^2 = b^2 = \beta w$ for all $z \in Z_0$. Hence the algebras with $\beta = 0$ and $\beta = 1$ are non-isomorphic (for $\beta = 0$ there are $16$ elements with zero square in $J$, while for $\beta = 1$ there are only $8$). Thus Orbits 1 and 2 give $2$ isomorphism classes each ($R_{3.1}, R_{3.2}$ and $R_{3.3}, R_{3.4}$).
\end{itemize}
In total, in Family 3 we obtain $2 + 2 + 1 + 1 + 1 = 7$ isomorphism classes. $\square$

\begin{table}[htbp]
\centering
\caption{Classification of algebras of Family 2 ($J^3 = 0$) and Family 3 ($J^3 \neq 0$) for $\mathrm{char}(R) = 2$, $\dim_{\mathbb{F}_2} W_0 = 1$}
\vspace{0.5em}
\small
\begin{tabular}{ |c|c|c|c|}
\hline
Class & Form $P(X, Y)$ & Products $a^2, ab, ba, b^2$ & $|\{x \in J : x^2 = 0\}|$ \\
\hline
$R_{2.1}$ & $0$ ($J^3 = 0$) & $a^2 = 0,\ ab = w,\ ba = 0,\ b^2 = 0$ & $24$ \\
\hline
$R_{2.2}$ & $0$ ($J^3 = 0$) & $a^2 = w,\ ab = w,\ ba = 0,\ b^2 = w$ & $8$ \\
\hline
$R_{3.1}$ & $X^3$ ($\beta = 0$) & $a^2 = c_1,\ ab = w,\ ba = 0,\ b^2 = 0$ & $16$ \\
\hline
$R_{3.2}$ & $X^3$ ($\beta = 1$) & $a^2 = c_1,\ ab = w,\ ba = 0,\ b^2 = w$ & $8$ \\
\hline
$R_{3.3}$ & $X^2 Y$ ($\beta = 0$) & $a^2 = c_2,\ ab = c_1 + w,\ ba = c_1,\ b^2 = 0$ & $16$ \\
\hline
$R_{3.4}$ & $X^2 Y$ ($\beta = 1$) & $a^2 = c_2,\ ab = c_1 + w,\ ba = c_1,\ b^2 = w$ & $8$ \\
\hline
$R_{3.5}$ & $X^2 Y + X Y^2$ & $a^2 = c_2,\ ab = c_1 + c_2 + w,\ ba = c_1 + c_2,\ b^2 = c_1$ & $8$ \\
\hline
$R_{3.6}$ & $X^3 + X^2 Y + X Y^2$ & $a^2 = c_1 + c_2,\ ab = c_1 + c_2 + w,\ ba = c_1 + c_2,\ b^2 = c_1$ & $8$ \\
\hline
$R_{3.7}$ & $X^3 + X^2 Y + Y^3$ & $a^2 = c_1 + c_2,\ ab = c_1 + w,\ ba = c_1,\ b^2 = c_2$ & $8$ \\
\hline
\end{tabular}
\end{table}

\section{Classification of Rings of Characteristic 4}

Let now $\mathrm{char}(R) = 4$. Then the element $e_0 = 2\cdot 1 = 1 + 1$ is a nonzero element of the central ideal $Z_0 \cong \mathbb{F}_2^3$. Any ring isomorphism $\varphi : R \xrightarrow{\sim} R'$ maps the unity $1$ to the unity $1'$, and therefore must map $e_0 = 2\cdot 1$ to $e_0' = 2\cdot 1'$. Moreover, by the distributive law, for any element $r \in J$ we have $2r = (1 + 1)r = e_0 r \in JZ_0 \subseteq W_0$. According to Theorem 1.4, the identities $2Z_0 = 0$ and $2J^2 = 0$ always hold in the radical $J$. In particular, for any $r \in J$ and $z \in Z_0$, since $z \in Z(R)$ and $rz \in Z_0$, we have $rz + zr = 2(rz) = 0$ and $z^2 = 0$, whence
$$
(r + z)^2 = r^2 + rz + zr + z^2 = r^2. \eqno(4.1)
$$
Similarly, for any $a, b \in J$, in view of $ab - ba = [a, b] \in W_0$ and $2(ba) \in 2J^2 = \{0\}$, we have $ab + ba = [a, b] + 2(ba) = [a, b]$. This means that all commutator and quadratic relations inside the radical $J$ obey the same formulas as in characteristic $2$, and the ring $R$ is uniquely recovered from the structure of the ring without unity $J$ and the distinguished element $e_0 = 2\cdot 1 \in Z_0 \setminus \{0\}$ with the consistency condition $2r = e_0 r$ ($r \in J$).

Depending on whether the element $e_0 = 2\cdot 1$ belongs to the annihilator $W_0$ or lies in $Z_0 \setminus W_0$, two fundamentally different cases arise.

\textbf{Case $2\cdot 1 \in W_0 \setminus \{0\}$ (rings with additive group $\mathbb{Z}_4 \oplus \mathbb{Z}_2^4$)}

If $e_0 = 2\cdot 1 \in W_0 \setminus \{0\}$, then for every $r \in J$ we have $2r = e_0 r = 0$ (since $W_0 \subseteq \mathrm{Ann}(J)$). Consequently, the ideal $J$ has exponent $2$ ($2J = 0$) and is a nilpotent $\mathbb{F}_2$-algebra from Families 1, 2, or 3 described in Section 3, and the ring $R$ is a central extension
$$
R = (\mathbb{Z}_4 \cdot 1 \oplus J) / \langle 2\cdot 1 - e_0 \rangle,
$$
whose additive group is isomorphic to $\mathbb{Z}_4 \oplus \mathbb{Z}_2^4$. Two such rings $(J, e_0)$ and $(J', e_0')$ are isomorphic if and only if there exists an isomorphism of $\mathbb{F}_2$-algebras $\varphi : J \xrightarrow{\sim} J'$ mapping $e_0$ to $e_0'$.

\textbf{4.1. Proposition.}
In the class of non-commutative centrally essential rings $R$ of order $64$ of characteristic $4$ with the condition $2\cdot 1 \in W_0 \setminus \{0\}$ there exist exactly $35$ pairwise non-isomorphic rings:

\textbf{1.} For $\dim_{\mathbb{F}_2} W_0 = 1$ ($W_0 = \{0, w\}$, $2\cdot 1 = w$) there are exactly $9$ isomorphism classes: $2$ classes with $J^3 = 0$ (corresponding to $R_{2.1}$ and $R_{2.2}$) and $7$ classes with $J^3 \neq 0$ (corresponding to $R_{3.1}$--$R_{3.7}$).

\textbf{2.} For $\dim_{\mathbb{F}_2} W_0 = 2$ ($W_0 = \langle s_1, s_2 \rangle_{\mathbb{F}_2}$) and $2\cdot 1 = s_1 \in [J, J]$ there are exactly $10$ isomorphism classes, corresponding to the classes $R_{1.1}$--$R_{1.10}$ of Table 1.

\textbf{3.} For $\dim_{\mathbb{F}_2} W_0 = 2$ and $2\cdot 1 \notin [J, J]$ (after choosing a basis, $2\cdot 1 = s_2$) there are exactly $16$ isomorphism classes, in one-to-one correspondence with all $16$ pairs $(\alpha, \beta) \in W_0 \times W_0$ in relations (3.1).

\textbf{Proof.}
\textbf{1.} If $\dim_{\mathbb{F}_2} W_0 = 1$, then the subgroup $W_0 = \{0, w\}$ contains a unique nonzero element $w$, which is fixed under any automorphism of the algebra $J$. Therefore fixing $2\cdot 1 = w$ gives a one-to-one correspondence with the isomorphism classes of algebras $J$ from Propositions 3.4 and 3.6, that is, $2 + 7 = 9$ classes.

\textbf{2.} Let $\dim_{\mathbb{F}_2} W_0 = 2$. In any algebra $J$ of Family 1 the commutator $[J, J] = \{0, s_1\}$ is a characteristic subspace. If $e_0 = 2\cdot 1 = s_1$, then the condition $\varphi(e_0) = e_0'$ holds automatically for any isomorphism $\varphi : J \xrightarrow{\sim} J'$, whence we obtain exactly $10$ isomorphism classes from Proposition 3.2.

\textbf{3.} If $e_0 = 2\cdot 1 \in W_0 \setminus [J, J] = \{s_2, s_1 + s_2\}$, then in the proof of Proposition 3.2 the only nontrivial change of the canonical basis ($\overline{b}' = \overline{a} + \overline{b}$) maps the basis vector $s_2$ to $s_2' = s_1 + s_2 \neq s_2$. Consequently, the requirement $s_2 = 2\cdot 1$ uniquely fixes the choice of both basis vectors $s_1 = [a, b]$ and $s_2 = 2\cdot 1$ in $W_0$, and hence also the choice of the cosets $a + Z_0 = \overline{L}_c^{-1}(s_1)$ and $b + Z_0 = \overline{L}_c^{-1}(s_2)$ in $J/Z_0$. By the invariance of the squares $(a + z)^2 = a^2 = \alpha$ and $(b + z)^2 = b^2 = \beta$, all $16$ pairs $(\alpha, \beta) \in W_0 \times W_0$ define pairwise non-isomorphic rings. $\square$

\textbf{Case $2\cdot 1 \in Z_0 \setminus W_0$ and $\dim_{\mathbb{F}_2} W_0 = 2$ (free $\mathbb{Z}_4$-algebras of rank 3)}

Let $\mathrm{char}(R) = 4$, $\dim_{\mathbb{F}_2} W_0 = 2$, and the element $c = 2\cdot 1$ lie in $Z_0 \setminus W_0$. Then $Z_0 = \langle 2\cdot 1, s_1, s_2 \rangle_{\mathbb{F}_2}$, and the right multiplication operator $L_c : J \to W_0$ coincides with the multiplication-by-$2$ operator: $L_c(r) = r(2\cdot 1) = 2r$ ($r \in J$). Choosing elements $a, b \in J \setminus Z_0$ such that $2a = ac = s_1$ and $2b = bc = s_2$, we see that the elements $1, a, b$ have additive order $4$ and generate $R$ as a free module of rank $3$ over the residue ring $\mathbb{Z}_4$:
$$
R = \mathbb{Z}_4 \cdot 1 \oplus \mathbb{Z}_4 a \oplus \mathbb{Z}_4 b \cong \mathbb{Z}_4^3, \,\, Z(R) = \mathbb{Z}_4 \cdot 1 \oplus 2\mathbb{Z}_4 a \oplus 2\mathbb{Z}_4 b, \,\, W_0 = \langle 2a, 2b \rangle.
$$
\textbf{4.2. Proposition.}
Every non-commutative centrally essential ring $R$ of order $64$ with the conditions $\mathrm{char}(R) = 4$, $\dim_{\mathbb{F}_2} W_0 = 2$, and $2\cdot 1 \in Z_0 \setminus W_0$ is a free $\mathbb{Z}_4$-algebra of rank $3$ with basis $\{1, a, b\}$ over $\mathbb{Z}_4$ and defining relations
$$
ab = 2a, \,\, ba = 0, \,\, a^2 = \alpha, \,\, b^2 = \beta, \,\, \text{where } \alpha, \beta \in W_0 = \{0, 2a, 2b, 2a + 2b\}. \eqno (4.2)
$$
There exist exactly $10$ pairwise non-isomorphic rings of this type, which are in one-to-one correspondence with the $10$ classes of Table 1 (under the identification $s_1 = 2a$, $s_2 = 2b$).

\textbf{Proof.} By Theorem 1.4 (item 3), for $\dim_{\mathbb{F}_2} W_0 = 2$ we always have $J^2 \subseteq W_0$. Setting $c = 2\cdot 1 \in Z_0 \setminus W_0$, note that in the proof of Proposition 3.2 the element $c$ was not subjected to any changes: the choice of $a, b \in J \setminus Z_0$ with $2a = ac = s_1 = [a, b]$ and $2b = bc = s_2$ and the annihilation of the product $ba = 0$ by the shift $a' = a + \xi(2\cdot 1)$, $b' = b + \eta(2\cdot 1)$ occur with $c = 2\cdot 1$ fixed. By identity (4.1), the squares $(a')^2 = a^2 = \alpha$ and $(b')^2 = b^2 = \beta$ do not depend on the choice of representatives in the cosets modulo $Z_0 = 2R$, and the change $\overline{b}' = \overline{a} + \overline{b}$ induces on the pairs $(\alpha, \beta) \in W_0^2$ the same involution $\tau$ as in Proposition 3.2. Consequently, we obtain the same $10$ isomorphism classes. $\square$

\textbf{Case $2\cdot 1 \in Z_0 \setminus W_0$ and $\dim_{\mathbb{F}_2} W_0 = 1$ (rings with group $\mathbb{Z}_4^2 \oplus \mathbb{Z}_2^2$)}

Let, finally, $\mathrm{char}(R) = 4$, $\dim_{\mathbb{F}_2} W_0 = 1$ ($W_0 = \{0, w\}$), and the canonical element $c_1 = 2\cdot 1$ belong to $Z_0 \setminus W_0$. Complement the vector $c_1 + W_0$ to a basis $\{c_1 + W_0, c_2 + W_0\}$ of the space $Z_0/W_0$ and choose the dual basis $\{a + Z_0, b + Z_0\}$ in $J/Z_0$:
$$
2a = ac_1 = w, \quad ac_2 = 0, \quad 2b = bc_1 = 0, \quad bc_2 = w.
$$
Then the elements $1$ and $a$ have additive order $4$ ($2\cdot 1 = c_1$, $2a = w$), while the elements $b$ and $c_2$ have additive order $2$, so that the additive group of the ring $R$ has the form:
$$
R = \mathbb{Z}_4 \cdot 1 \oplus \mathbb{Z}_4 a \oplus \mathbb{Z}_2 b \oplus \mathbb{Z}_2 c_2 \cong \mathbb{Z}_4 \oplus \mathbb{Z}_4 \oplus \mathbb{Z}_2 \oplus \mathbb{Z}_2.
$$
By Lemma 3.5 the multiplication table of the generators $a$ and $b$ is given by relations (3.4):
$$
a^2 = \lambda_1(2\cdot 1) + \lambda_2 c_2 + \alpha(2a), \,\, ab = \lambda_2(2\cdot 1) + \lambda_3 c_2 + 2a, \eqno (4.3)
$$
$$
ba = \lambda_2(2\cdot 1) + \lambda_3 c_2, \,\, b^2 = \lambda_3(2\cdot 1) + \lambda_4 c_2 + \beta(2a),
$$
where $\lambda_1, \lambda_2, \lambda_3, \lambda_4, \alpha, \beta \in \mathbb{F}_2$.

\textbf{4.3. Proposition.}
In the class of non-commutative centrally essential rings $R$ of order $64$ with the conditions $\mathrm{char}(R) = 4$, $\dim_{\mathbb{F}_2} W_0 = 1$, and $2\cdot 1 \in Z_0 \setminus W_0$ there exist exactly $24$ pairwise non-isomorphic rings:

\textbf{1.} For $J^2 \subseteq W_0$ ($(\lambda_1, \lambda_2, \lambda_3, \lambda_4) = (0, 0, 0, 0)$, the nilpotency index of $J$ equals $3$) there are exactly $3$ isomorphism classes, corresponding to the pairs $(\alpha, \beta) \in \{(0, 0), (0, 1), (1, 1)\}$.

\textbf{2.} For $J^2 \nsubseteq W_0$ ($(\lambda_1, \lambda_2, \lambda_3, \lambda_4) \neq (0, 0, 0, 0)$, the nilpotency index of $J$ equals $4$) there are exactly $21$ isomorphism classes, the complete list of which is given in Table 3.

\textbf{Proof.} Since the elements $c_1 = 2\cdot 1 \in Z_0 \setminus W_0$ and $w \in W_0 \setminus \{0\}$ are uniquely determined by the ring structure of $R$, any ring isomorphism of this class must leave $c_1$ and $w$ fixed: $c_1' = c_1$ and $w' = w$. Consequently, an arbitrary admissible change of the element $c_2 \in Z_0 \setminus \langle c_1, w \rangle$ has the form
$$
c_2' = c_2 + \lambda c_1 + \mu_2 w \qquad (\lambda, \mu_2 \in \mathbb{F}_2).
$$
The duality condition of the basis $\{a' + Z_0, b' + Z_0\}$ to the space $\{c_1 + W_0, c_2' + W_0\}$ ($a'c_1 = w$, $a'c_2' = 0$, $b'c_1 = 0$, $b'c_2' = w$) uniquely determines the transformation in $J/Z_0$:
$$
a' \equiv a + \lambda b \pmod{Z_0}, \qquad b' \equiv b \pmod{Z_0}.
$$
Thus on the set of parameters $(\lambda_1, \lambda_2, \lambda_3, \lambda_4; \alpha, \beta) \in \mathbb{F}_2^4 \times \mathbb{F}_2^2$ there acts the group of admissible basis changes $G = \{\mathrm{id}, \sigma, T_1, \sigma T_1\} \cong C_2 \times C_2$ of order $4$, generated by two commuting involutions:

\textbf{Shift $T_1$} ($\lambda = 0, \mu_2 = 1$): $c_1' = c_1$, $c_2' = c_2 + w$, $a' = a + \lambda_3 c_2$, $b' = b$. Substituting $c_2 = c_2' + w$ into (4.3) preserves the quadruple $(\lambda_1, \lambda_2, \lambda_3, \lambda_4)$ and transforms $(\alpha, \beta)$ by the rule $T_1(\alpha, \beta) = (\alpha + \lambda_2, \, \beta + \lambda_4)$.

\textbf{Borel transformation $\sigma$} ($\lambda = 1, \mu_2 = 0$): $c_1' = c_1$, $c_2' = c_2 + c_1$, $a' = a + b + \beta c_2'$, $b' = b$. Substituting $c_2 = c_2' + c_1'$ into the expressions for $(b')^2 = b^2$, $b'a' = ba + b^2 + \beta w$, and $(a')^2 = a^2 + b^2 + w$, we find:
$$
\sigma(\lambda_1, \lambda_2, \lambda_3, \lambda_4) = (\lambda_1 + \lambda_2 + \lambda_3 + \lambda_4, \, \lambda_2 + \lambda_4, \, \lambda_3 + \lambda_4, \, \lambda_4),
$$
$$
\sigma(\alpha, \beta) = (\alpha + \beta + 1, \, \beta).
$$
\textbf{1.} Let $J^2 \subseteq W_0$, that is, $(\lambda_1, \lambda_2, \lambda_3, \lambda_4) = (0, 0, 0, 0)$. Then $T_1$ acts on $(\alpha, \beta)$ identically, and the involution $\sigma(\alpha, \beta) = (\alpha + \beta + 1, \beta)$ interchanges the pairs $(0, 0) \leftrightarrow (1, 0)$ and leaves the pairs $(0, 1)$ and $(1, 1)$ fixed. This gives exactly $3$ isomorphism classes.

\textbf{2.} Let $J^2 \nsubseteq W_0$, that is, $(\lambda_1, \lambda_2, \lambda_3, \lambda_4) \in \mathbb{F}_2^4 \setminus \{(0, 0, 0, 0)\}$. Let us compute the number of orbits of the group $G = \{\mathrm{id}, \sigma, T_1, \sigma T_1\}$ on the set $X = (\mathbb{F}_2^4 \setminus \{0\}) \times \mathbb{F}_2^2$ of $15 \times 4 = 60$ elements by Burnside's lemma. We have $|\mathrm{Fix}(\mathrm{id})| = 60$. Further, the element $T_1$ fixes a set $(\lambda_1, \lambda_2, \lambda_3, \lambda_4; \alpha, \beta)$ if and only if $\lambda_2 = \lambda_4 = 0$. This condition is satisfied by $2^2 - 1 = 3$ nonzero vectors $(\lambda_1, 0, \lambda_3, 0)$ and all $4$ pairs $(\alpha, \beta)$, whence $|\mathrm{Fix}(T_1)| = 3 \times 4 = 12$.

The involution $\sigma$ fixes a vector $(\lambda_1, \lambda_2, \lambda_3, \lambda_4)$ if and only if $\lambda_4 = 0$ and $\lambda_2 = \lambda_3$ ($3$ nonzero vectors: $(1, 0, 0, 0)$, $(0, 1, 1, 0)$, $(1, 1, 1, 0)$), and fixes a pair $(\alpha, \beta)$ if and only if $\beta = 1$ ($2$ pairs: $(0, 1)$ and $(1, 1)$). Hence $|\mathrm{Fix}(\sigma)| = 3 \times 2 = 6$.

The involution $\sigma T_1$ fixes the same $3$ vectors with $\lambda_4 = 0$, $\lambda_2 = \lambda_3$, and acts on the pairs $(\alpha, \beta)$ as $(\alpha, \beta) \mapsto (\alpha + \beta + 1 + \lambda_2, \beta)$, having exactly $2$ fixed pairs ($\beta = 1 + \lambda_2$) for each of the $3$ vectors. Hence $|\mathrm{Fix}(\sigma T_1)| = 3 \times 2 = 6$.

By Burnside's formula the number of orbits (isomorphism classes) for $J^2 \nsubseteq W_0$ equals:
$$
N = \frac{1}{|G|}\bigl(|\mathrm{Fix}(\mathrm{id})| + |\mathrm{Fix}(T_1)| + |\mathrm{Fix}(\sigma)| + |\mathrm{Fix}(\sigma T_1)|\bigr) = \frac{60 + 12 + 6 + 6}{4} = 21.
$$
Explicit representatives of all $9$ orbits of vectors $(\lambda_1, \lambda_2, \lambda_3, \lambda_4)$ under $\sigma$ and the corresponding admissible pairs $(\alpha, \beta)$ are given in Table 3.
$\square$

\begin{table}[htbp]
\centering
\caption{Classification of rings for $\operatorname{char}(R) = 4$, $\dim_{\mathbb{F}_2} W_0 = 1$, $2\cdot 1 = c_1 \in Z_0 \setminus W_0$ ($24$ classes)}
\vspace{0.5em}
\small
\begin{tabular}{ |c|c|c|c|c|}
\hline
Subfamily & Orbit $(\lambda_1, \lambda_2, \lambda_3, \lambda_4)$ w.r.t. $\sigma$ & Representative & Classes $(\alpha, \beta)$ & Number \\
\hline
$J^3 = 0$ & $\{(0,0,0,0)\}$ & $(0,0,0,0)$ & $(0,0), (0,1), (1,1)$ & $3$ \\
\hline
$J^3 \neq 0$ (orb. 1) & $\{(1,0,0,0)\}$ & $(1,0,0,0)$ & $(0,0), (0,1), (1,1)$ & $3$ \\
\hline
$J^3 \neq 0$ (orb. 2) & $\{(0,1,1,0)\}$ & $(0,1,1,0)$ & $(0,0), (0,1)$ & $2$ \\
\hline
$J^3 \neq 0$ (orb. 3) & $\{(1,1,1,0)\}$ & $(1,1,1,0)$ & $(0,0), (0,1)$ & $2$ \\
\hline
$J^3 \neq 0$ (orb. 4) & $\{(0,1,0,0), (1,1,0,0)\}$ & $(0,1,0,0)$ & $(0,0), (0,1)$ & $2$ \\
\hline
$J^3 \neq 0$ (orb. 5) & $\{(0,0,1,0), (1,0,1,0)\}$ & $(0,0,1,0)$ & $(0,0), (1,0), (0,1), (1,1)$ & $4$ \\
\hline
$J^3 \neq 0$ (orb. 6) & $\{(0,0,0,1), (1,1,1,1)\}$ & $(0,0,0,1)$ & $(0,0), (1,0)$ & $2$ \\
\hline
$J^3 \neq 0$ (orb. 7) & $\{(1,0,0,1), (0,1,1,1)\}$ & $(1,0,0,1)$ & $(0,0), (1,0)$ & $2$ \\
\hline
$J^3 \neq 0$ (orb. 8) & $\{(0,1,0,1), (0,0,1,1)\}$ & $(0,1,0,1)$ & $(0,0), (0,1)$ & $2$ \\
\hline
$J^3 \neq 0$ (orb. 9) & $\{(1,0,1,1), (1,1,0,1)\}$ & $(1,0,1,1)$ & $(0,0), (1,0)$ & $2$ \\
\hline
\textbf{Total} & & & & $\mathbf{24}$ \\
\hline
\end{tabular}
\end{table}

\end{document}